\documentclass[a4paper, 11pt]{amsproc}
\usepackage{amsmath}    
\usepackage{amsthm}     
\usepackage{amssymb}    
\usepackage{mathtools}  
\usepackage{amstext}    

\usepackage{mathrsfs}   
\usepackage{bm}         
\usepackage{shuffle}    
\usepackage[all]{xy}    
\usepackage[T1]{fontenc}
\usepackage{lmodern}

\usepackage{amsrefs}    
\usepackage{ascmac}     
\usepackage{comment}    
\usepackage{ytableau}   
\usepackage{here}       
\usepackage{fullpage}   
\usepackage{xcolor}     

\usepackage[utf8]{inputenc} 

\usepackage{hyperref}   
\hypersetup{
 colorlinks=true,       
 linkcolor=blue,        
 citecolor=blue,        
 filecolor=blue,        
 urlcolor=blue          
}

\usepackage[capitalize,nameinlink,noabbrev,nosort]{cleveref}

\DeclareFontEncoding{OT2}{}{}
\DeclareFontSubstitution{OT2}{cmr}{m}{l}
\DeclareFontFamily{OT2}{cmr}{\hyphenchar\font45 }
\DeclareFontShape{OT2}{cmr}{m}{l}{%
  <5><6><7><8><9>gen*wncyr%
  <10><10.95><12><14.4><17.28><20.74><24.88>wncyr10}{}
\DeclareMathAlphabet{\mathcyr}{OT2}{cmr}{m}{l}
\DeclareMathAlphabet{\mathcyb}{OT2}{cmr}{b}{l}
\SetMathAlphabet{\mathcyr}{bold}{OT2}{cmr}{b}{l}

\makeatletter
\@namedef{subjclassname@2020}{%
  \textup{2020} Mathematics Subject Classification}
\makeatother

\newtheorem{theoremcounter}{Theorem Counter}[section]
\theoremstyle{definition}
\newtheorem{definition}[theoremcounter]{Definition}
\newtheorem{remark}[theoremcounter]{Remark}
\newtheorem{example}[theoremcounter]{Example}
\theoremstyle{plain}
\newtheorem{lemma}[theoremcounter]{Lemma}

\newtheorem{corollary}[theoremcounter]{Corollary}

\newtheorem{theorem}[theoremcounter]{Theorem}

\theoremstyle{definition}

\theoremstyle{plain}

\numberwithin{equation}{section}

\allowdisplaybreaks[1]
\usepackage{enumitem}

\begin{document}

\title{On the uniform distribution modulo $1$ of zeros and $a$-points of zeta functions}

\author{Hideki Murahara}
\address[Hideki Murahara]{The University of Kitakyushu,  4-2-1 Kitagata, Kokuraminami-ku, Kitakyushu, Fukuoka, 802-8577, Japan}
\email{hmurahara@mathformula.page}

\author{Tomokazu Onozuka}
\address[Tomokazu Onozuka]{Division of Mathematical Sciences, Department of Integrated Science and Technology, Faculty of Science and Technology, Oita University, 700 Dannoharu, Oita, 870-1192, Japan} 
\email{t-onozuka@oita-u.ac.jp}


\begin{abstract}
Fujii gave five sufficient conditions for the uniform distribution modulo $1$ of the sequence $(u f(\gamma_n))$, where $\gamma_n$ runs over the imaginary parts of the non-trivial zeros of the Riemann zeta function. 
In this paper, we provide four sufficient conditions for the uniform distribution modulo $1$ that apply to a much broader class of sequences.
Our method relies solely on the asymptotic behavior of the counting function and does not require the intricate calculations concerning the Riemann zeta function employed by Fujii in his paper.
As applications, we prove the uniform distribution modulo $1$ of $(u f(x_n))$ for various sequences $(x_n)$, including the non-trivial $a$-points of the derivatives of the Riemann zeta function, the non-trivial zeros of the derivatives of Dirichlet $L$-functions, and the non-trivial zeros of functions in the Selberg class. 
Furthermore, by applying the Erd\H{o}s--Tur\'an inequality, we obtain an upper bound for the discrepancy of these sequences.
\end{abstract}


\maketitle

\section{Introduction}
For $\Re(s)>1$, the Riemann zeta function $\zeta(s)$ is defined by the following Dirichlet series:
\[
\zeta(s)=\sum_{n=1}^\infty\frac{1}{n^s}.
\]
This function can be meromorphically continued to the entire complex plane. 
It is well known that $\zeta(s)$ is intimately related to prime numbers, and investigating the distribution of its zeros provides crucial information about the distribution of primes. 
The zeros of the Riemann zeta function are broadly classified into the well-known trivial zeros and the elusive non-trivial zeros. The non-trivial zeros are known to lie in the critical strip, where the real part is strictly between $0$ and $1$, and they are distributed symmetrically with respect to the critical line $\Re(s)=1/2$. 
The Riemann Hypothesis, one of the most famous unsolved problems in mathematics, asserts that all non-trivial zeros lie exactly on this critical line.

Therefore, we are particularly interested in the behavior of the Riemann zeta function in the right half of the critical strip, namely $(1/2, 1)$. 
However, $\zeta(s)$ exhibits highly elusive and complex behavior in this region, both in its value distribution and in the distribution of its zeros. 
A striking manifestation of the former is Voronin's universality theorem.
\begin{theorem}[Voronin's Universality Theorem \cite{Voronin1975}]
Let $K$ be a compact subset of the strip $1/2 < \Re(s) < 1$ with connected complement. Let $f(s)$ be a continuous function on $K$ which is analytic in the interior of $K$ and has no zeros in $K$. 
Then, for any $\epsilon>0$, we have
\[ 
\liminf_{T \to \infty} \frac{1}{T} \operatorname{meas}
\left\{ \tau \in [0, T] ~\middle|~ \max_{s \in K} |\zeta(s+i\tau) - f(s)| < \epsilon \right\}>0, 
\]
where $\operatorname{meas}$ denotes the Lebesgue measure.
\end{theorem}
This theorem asserts that the Riemann zeta function can approximate any non-vanishing analytic function on compact subsets of the critical strip with positive lower density, indicating that its value distribution is far from simple.
The distribution of its zeros is equally elusive. 
Let $\{x\}=x - \lfloor x \rfloor$ denote the fractional part of a real number $x$, where $\lfloor x \rfloor$ is the greatest integer less than or equal to $x$. We then define uniform distribution modulo $1$ as follows: 
\begin{definition}[Definition 1.1 in Chapter 1 of \cite{KuipersNiederreiter1974}]
A sequence $(x_n)_{n=1}^\infty$ of real numbers is said to be uniformly distributed modulo $1$ if for all real numbers $a$ and $b$ with $0\le a<b\le 1$, we have
\[
\lim_{N\to\infty} 
\frac{\#\{1\le n\le N\mid a\le \{x_n\}<b\}}{N}
=b-a.
\]
\end{definition}

Let $\rho_n=\beta_n +i\gamma_n$ denote the non-trivial zeros of the Riemann zeta function with $\gamma_n>0$, counted with multiplicity and ordered so that $0 < \gamma_1 \le \gamma_2 \le \cdots$. 
In 1956, Rademacher~\cite{Rademacher1974} showed, under the assumption of the Riemann Hypothesis, that for any fixed non-zero real number $u$, the sequence $(u\gamma_n)$ is uniformly distributed modulo $1$. 
Unaware of Rademacher's work, Elliott~\cite{Elliott1972} independently proved  the same conditional result. 
His proof relied on a conditional estimate due to Selberg for the second moment of $S(t)$. 
In reality, however, an unconditional estimate—albeit with a slightly weaker error term—had already been obtained by Selberg~\cite{zbMATH03099795} himself. 
According to an ``Added note'' in Elliott's paper, this fact was pointed out to him by Selberg directly during a meeting in Oberwolfach, rendering his result unconditional.
Subsequently, an unconditional proof of this uniform distribution was given independently by Hlawka~\cite{HlawkaEdmund1976}. 

Akbary and Murty~\cite{AkbaryMurty2008} established, under the average density hypothesis, the uniform distribution modulo $1$ of the non-trivial zeros of functions in a class of Dirichlet series containing the Selberg class.
This line of research was naturally expanded to $a$-points. 
Recall that for a given complex number $a$, $a$-points of a function $F(s)$ are the roots of the equation $F(s)=a$. 
Steuding~\cite{Steuding2012,Steuding2014} proved the uniform distribution modulo $1$ for the non-trivial $a$-points of the Riemann zeta function.
Jakhlouti, Mazhouda, and Steuding~\cite{JakhloutiMazhoudaSteuding2015} generalized this to the $a$-points of Selberg class $L$-functions. 

Investigations into the derivatives of these functions soon followed. 
Lee, the second author, and Suriajaya~\cite{LeeOnozukaSuriajaya2016} obtained the corresponding result for the non-trivial $a$-points of the derivatives of the Riemann zeta function, which was recently extended by Mekkaoui and Mazhouda~\cite{MekkaouiMazhouda2024} to the $a$-points of the derivatives of functions in the Selberg class.

Beyond the Riemann zeta function and the Selberg class, 
similar distribution properties have been explored for other zeta functions. 
For instance, Schmeller~\cite{Schmeller2018} demonstrated the uniform distribution modulo $1$ for the ordinates of the non-trivial zeros of the Epstein zeta function under a certain assumption, while Garunk{\v{s}}tis and Panavas~\cite{GarunkstisPanavas2022} established it for the non-trivial zeros of the Lerch zeta function. 
Furthermore, Garunk{\v{s}}tis, Steuding, and {\v{S}}im{\.e}nas~\cite{GarunkstisSteuding2014} obtained similar results for the $a$-points of the Selberg zeta function.
These studies highlight an interesting contrast: although the imaginary parts of these zeros (such as $\gamma_n$) grow irregularly and have no simple pattern, under appropriate hypotheses, they are all uniformly distributed modulo $1$.

However, it is important to note that such an apparent irregularity does not necessarily imply a uniform distribution modulo $1$. 
For instance, Wintner~\cite{zbMATH02531574} proved that the sequence $(\log p_n)$ is not uniformly distributed modulo $1$, where $(p_n)$ is the sequence of prime numbers in ascending order. Motivated by this result, Rehberg~\cite{zbMATH07099390} investigated what kinds of functions $f$ would make the sequence modulo $1$ uniformly distributed. 
To be precise, Rehberg provided conditions on a function $f$ such that the sequence $(uf(q_n))$ is uniformly distributed modulo $1$ for any non-zero real number $u$. Here, $(q_n)$ is a sequence of real numbers whose counting function satisfies
\[
\sum_{q_n \le x} 1
=\alpha \mathrm{Li}(x) +O\left( \frac{x}{(\log x)^k} \right)
\]
for some constant $\alpha>0$ and any fixed $k>1$. 
In this formula, $\mathrm{Li}(x)$ denotes the logarithmic integral function.

On the other hand, Fujii~\cite{AkioFujii1982distribution2} provided conditions on a function $f$ such that the real sequence $(uf(\gamma_n))$ is uniformly distributed modulo $1$. From Fujii's results, it follows that, unlike the case of prime numbers, uniform distribution holds even when $f(x)=\log x$. That is, he showed that the sequence $(\log \gamma_n)$ is uniformly distributed modulo $1$.
To state further examples derived from his work, we define $\log_k x$ to be the $k$-fold iterated logarithm, defined recursively by $\log_1 x=\log x$ and $\log_k x=\log(\log_{k-1} x)$ for $k \ge 2$. As concrete consequences of his results, we have the following:
\begin{theorem}[Fujii~\cite{AkioFujii1982distribution2}]
Sequences such as
$(\gamma_n \log \gamma_n/\log_k \gamma_n)$,
$(\gamma_n (\log \gamma_n)^b)$,
$(\gamma_n^{b'})$,
$((\log \gamma_n)^{b''})$,
and $(\log \gamma_n \cdot \log_k \gamma_n)$
are uniformly distributed modulo $1$, where $k$ is a positive integer, $b<1$, $0<b'\le1$, and $b''>1$.
\end{theorem}

Regarding the method of proof, Rehberg's results used the concept of discrepancy. The discrepancy $D_N$ is defined as follows:
\begin{definition}[Definition 1.1 in Chapter 2 of \cite{KuipersNiederreiter1974}]
For a sequence $(x_n)_{n=1}^N$ of real numbers, the number
\[
D_N=D_N(x_1,x_2,\ldots,x_N)
=\sup_{0\le a<b\le1} 
\left| 
\frac{\#\{1\le n\le N\mid a\le \{x_n\}<b\}}{N} - (b-a)
\right|
\]
is called the discrepancy of $(x_n)_{n=1}^N$.
\end{definition}
By \cite[Chapter 2, Theorem 1.1]{KuipersNiederreiter1974}, it is known that a sequence of real numbers $(x_n)_{n=1}^\infty$ is uniformly distributed modulo $1$ if and only if $\lim_{N \to \infty} D_N=0$. 
To estimate this discrepancy and effectively apply this equivalence, Rehberg also employed the Erd\H{o}s--Tur\'an inequality. 
\begin{theorem}[Erd\H{o}s--Tur\'an inequality, (2.42) in Chapter 2 of \cite{KuipersNiederreiter1974}] \label{thm:erdos_turan}
Let $x_1, x_2, \ldots, x_N$ be real numbers, and let $D_N$ denote the discrepancy of the sequence $(x_n)_{n=1}^N$. 
Then, for any positive integer $H$, there exists an absolute constant $C>0$ such that
\[
D_N \le C 
\left( 
 \frac{1}{H} 
 +\sum_{h=1}^H \frac{1}{h} 
 \left| \frac{1}{N} \sum_{n=1}^N e^{2\pi i h x_n} \right| 
\right).
\]
\end{theorem}
By evaluating the exponential sums on the right-hand side of this inequality, Rehberg derived four conditions on the function $f$ to ensure that the sequence $(uf(q_n))$ is uniformly distributed modulo $1$.

In contrast, Fujii gave a proof using Weyl's criterion. 
\begin{lemma}[Weyl's criterion, Theorem 2.1 in Chapter 1 of \cite{KuipersNiederreiter1974}] \label{lem:Weyl}
Let $(x_k)_{k=1}^\infty$ be a sequence of real numbers. 
Then the sequence $(x_k)$ is uniformly distributed modulo $1$ if and only if
\[
\lim_{N \to \infty} 
\frac{1}{N} 
\sum_{k=1}^N e^{2\pi i h x_k} 
=0 
\quad \text{for all } h\in\mathbb{Z} \setminus \{0\}.
\]
\end{lemma}
Applying this criterion, Fujii established five conditions on the function $f$ for the real sequence $(uf(\gamma_n))$ to be uniformly distributed modulo $1$.

Regarding the quantitative aspect, the discrepancy of the simpler sequence $(u\gamma_n)$ has been investigated. 
Let $D_{N(T)}$ denote the discrepancy of the sequence $(u\gamma_n)_{\gamma_n \le T}$ modulo $1$. 
Previously, Hlawka~\cite{HlawkaEdmund1976} proved that $D_{N(T)}=O(1/\log T)$ under the assumption of the Riemann Hypothesis. 
Fujii~\cite{AkioFujii1976no3} unconditionally established the bound $D_{N(T)}=O(\log\log T / \log T)$ for each non-zero real number $u$. 
Further deep investigations into the discrepancy of the non-trivial zeros were conducted by Ford, Soundararajan, and Zaharescu~\cite{FordSoundararajanZaharescu2009}.

In this paper, we give a general criterion for the uniform distribution modulo $1$ of sequences of real numbers. 
Our method relies on the asymptotic behavior of their counting functions rather than their specific arithmetic properties. 
Specifically, we focus on sequences $(x_n)_{n=1}^\infty$ tending to infinity whose counting function $N(T,(x_n)_{n=1}^\infty)\coloneqq\sum_{x_n \le T} 1$ satisfies the asymptotic formula
\[
N(T,(x_n)_{n=1}^\infty)=C_1 T (\log T-\log C_2) +O(\log T)
\]
as $T \to \infty$, where $C_1>0$ and $C_2>0$ are constants. 
As we will see, this broad framework encompasses many sequences in analytic number theory, including the $a$-points of the derivatives of the Riemann zeta function, the zeros of the derivatives of Dirichlet $L$-functions, and the zeros of functions in the Selberg class.

Previously, Fujii established a set of five conditions for the uniform distribution modulo $1$ of the sequence $(u f(\gamma_n))$. 
His proof, however, relies heavily on extensive computations on the critical line that are specific to the Riemann zeta function. 
Since performing analogous computations for general sequences (such as the $a$-points of derivatives) is difficult, his original conditions cannot be directly adapted to our broader context. 
To overcome this obstacle, we introduce four accessible conditions. 
Without relying on specific computations on the critical line, our method provides a general criterion applicable to a wide variety of sequences.

Before stating our main result, we introduce some notation and conventions. 
We use the standard Landau notation: for functions $f$ and $g$ with $g(x)\neq 0$ for sufficiently large $x$, we write $f(x)=O(g(x))$ as $x\to\infty$ if there exists a constant $C>0$ such that $|f(x)|\le C|g(x)|$ for all sufficiently large $x$, and $f(x)=o(g(x))$ as $x\to\infty$ if $\lim_{x\to\infty} f(x)/g(x)=0$. 
Hereinafter, all asymptotic statements and conditions are to be understood as $t \to \infty$ (or $T \to \infty$), unless otherwise specified.

Let $c$ be a sufficiently large constant. We introduce the following four conditions on a function $f$: 
\begin{enumerate}[label=(C\arabic*)]
\item $f(t)$ is a continuously differentiable real-valued function on $(c,\infty)$.
\item $\frac{f'(t)}{\log (et/C_2)}$ is monotonic on $[c,\infty)$ and remains strictly positive or strictly negative.
\item $f'(T)^{-1}=o(T)$.
\item $f(T)=o(T)$.
\end{enumerate}

We can now state our main criterion.
\begin{theorem}\label{main}
Let $C_1>0$ and $C_2>0$ be constants. 
Let $(x_n)_{n=1}^\infty$ be a sequence of real numbers whose counting function  $N(T,(x_n)_{n=1}^\infty)=\sum_{x_n \le T} 1$ satisfies
\[
N(T,(x_n)_{n=1}^\infty)=C_1 T (\log T-\log C_2) +O(\log T).
\]
Let $u$ be any non-zero real number. If a function $f$ satisfies the four conditions \textup{(C1)}--\textup{(C4)} above, then the sequence $(u f(x_n))_{x_n>c}$ is uniformly distributed modulo $1$.
\end{theorem}

Next, we apply this general framework to the specific case of the $k$-th derivative $\zeta^{(k)}(s)$. 
Let $\rho_{n}^{(k)}=\beta_{n}^{(k)} +i\gamma_{n}^{(k)}$ denote its non-trivial zeros, counted with multiplicity and ordered such that $0 < \gamma_{1}^{(k)} \le \gamma_{2}^{(k)} \le \cdots$. 
It is known~\cite{Akatsuka2012, Suriajaya2015} that, assuming the Riemann Hypothesis, the error term in the corresponding zero-counting function 
$N^{(k)}_0(T)\coloneqq\sum_{c<\gamma_{n}^{(k)}<T}1$ 
admits a sharper bound. 
This improved estimate allows us to relax condition \textup{(C4)} to the following weaker requirement:
\begin{enumerate}[label=(C4)']
\item $f(T)=o(T \sqrt{\log\log T})$.
\end{enumerate}
We thus obtain the following conditional result, where the value $C_2=4\pi e$ corresponds to $n_{0,k}=2$ for $k\ge 1$ in the Riemann--von Mangoldt formula for $\zeta^{(k)}(s)$.
\begin{theorem}\label{main2}
Assume the Riemann Hypothesis. Let $u$ be any non-zero real number and $k$ a positive integer. If a function $f$ satisfies conditions \textup{(C1)}--\textup{(C3)} with $C_2=4\pi e$, together with \textup{(C4)'}, then the sequence $(uf(\gamma_{n}^{(k)}))_{\gamma_{n}^{(k)}>c}$ is uniformly distributed modulo $1$.
\end{theorem}

Finally, by adapting Fujii's quantitative arguments to our generalized framework and applying the Erd\H{o}s--Tur\'an inequality, we obtain the following upper bound for the discrepancy.
\begin{theorem} \label{thm:discrepancy}
Let $u$ be any non-zero real number. 
Assume that a sequence $(x_n)$ and a function $f$ satisfy the conditions of Theorem \ref{main}. 
Let $D_{N(T)-N(c)}$ denote the discrepancy of the sequence $(u f(x_n))_{c < x_n \le T}$ modulo $1$. 
Then we have
\begin{align*}
D_{N(T)-N(c)}
&=O\left(\frac{1}{T |f'(T)|}\right) 
+O\left(\frac{\log T}{T}\right) 
+O\biggl(\left( \frac{T }{|f(T)| } \right)^{1/2}\biggr).
\end{align*}
\end{theorem}
Our framework depends only on the asymptotic counting formula $C_1 T(\log T-\log C_2)+O(\log T)$, and applies in principle to any zeta or $L$-function whose zero-counting (or $a$-point counting) function admits this shape. 
Further questions concern unconditional improvements of the discrepancy bound and the extension of the present method to the Lerch and Epstein zeta functions.

\section{Examples}
In this section, we present explicit examples of sequences $(x_n)$ and functions $f$ that fit into our general framework. The results below are immediate consequences of our main theorem.
\subsection{Examples of sequences $(x_n)$}
As our primary application of Theorem \ref{main}, we consider the sequence of $a$-points of the $k$-th derivative of the Riemann zeta function:
\[
\zeta^{(k)}(s)=\sum_{n=1}^\infty\frac{(-\log n)^k}{n^s}.
\]
 Let $\gamma_{a,n}^{(k)}$ denote the imaginary parts of these $a$-points, ordered such that $1 \le \gamma_{a,1}^{(k)} \le \gamma_{a,2}^{(k)} \le \cdots$. 
Their counting function is known to satisfy the required asymptotic formula with $C_1=1/(2\pi)$ and $C_2=2\pi e n_{a,k}$, where the constant $n_{a,k}$ is explicitly given by
\[
n_{a,k}=\begin{cases} 2 & \text{if }(a,k)=(1,0)\text{ or }a=0, k \ge 1, \\ 1 & \text{otherwise.}\end{cases}
\]
These values have been established through a series of foundational works: the case $a=k=0$ is the classical Riemann-von Mangoldt formula; the cases for $k=0$ with $a \ne 0$ are due to Landau~\cite{BohrLandauLittlewood1913} in a joint paper with Bohr and Littlewood; the cases for $a=0$ with $k \ge 1$ were obtained by Berndt~\cite{Berndt1970}; and the remaining cases for $a \ne 0$ with $k \ge 1$ were established by the second author~\cite{Onozuka2017}. 
Thus, we immediately obtain the following result.  
\begin{corollary}\label{cor:zeta_derivative}
Let $u$ be any non-zero real number. If a function $f$ satisfies the conditions \textup{(C1)}--\textup{(C4)} with $C_2=2\pi e n_{a,k}$, then the sequence $(uf(\gamma_{a,n}^{(k)}))_{\gamma_{a,n}^{(k)}>c}$ is uniformly distributed modulo $1$.
\end{corollary}

Our general criterion (Theorem \ref{main}) can also be applied to other objects in analytic number theory. As another application, we consider the zeros of the derivatives of Dirichlet $L$-functions. Let $\chi$ be a primitive Dirichlet character modulo $q>2$. In the half-plane $\Re(s)>1$, the associated Dirichlet $L$-function is defined by the series $L(s, \chi)=\sum_{n=1}^\infty \chi(n) n^{-s}$, and for any positive integer $k$, its $k$-th derivative is given by
\[
L^{(k)}(s, \chi)=(-1)^k \sum_{n=1}^{\infty} \frac{\chi(n) (\log n)^k}{n^s}.
\]
Let $\gamma_{n}^{(k)}(\chi)$ denote the imaginary parts of the zeros of $L^{(k)}(s, \chi)$ in the region $-q^K < \Re(s) < \sigma_k$, where $K$ and $\sigma_k$ are positive constants related to its zero-free region. We order these zeros by their absolute distance from the real axis such that $0 < |\gamma_{1}^{(k)}(\chi)| \le |\gamma_{2}^{(k)}(\chi)| \le \cdots$. As established by Y{\i}ld{\i}r{\i}m~\cite{Yildirim1996}, the counting function $N_{\chi}^{(k)}(T)$ for the zeros satisfying $|\gamma_{n}^{(k)}(\chi)| \le T$ satisfies the asymptotic formula
\[
N_{\chi}^{(k)}(T)=\frac{T}{\pi} \log \left(\frac{qT}{2\pi em}\right) +O(\log T),
\]
where $m>1$ is the smallest prime not dividing $q$. 
By rewriting the logarithmic term as $\log T - \log(2\pi e m / q)$, we see that this precisely matches our required asymptotic form with $C_1=1/\pi$ and $C_2=2\pi e m / q$. Therefore, we have the following result:
\begin{corollary}\label{cor:dirichlet_derivative}
Let $u$ be any non-zero real number. If a function $f$ satisfies the conditions \textup{(C1)}--\textup{(C4)} with $C_2=2\pi e m / q$, then the sequence $(u f(|\gamma_{n}^{(k)}(\chi)|))_{|\gamma_{n}^{(k)}(\chi)|>c}$ is uniformly distributed modulo $1$.
\end{corollary}
As our final application to specific sequences, we apply our criterion to the non-trivial zeros of functions in the Selberg class $\mathcal{S}$. For a function $L \in \mathcal{S}$, let $\gamma_{L,n}$ denote the imaginary parts of its non-trivial zeros, ordered such that $0 < \gamma_{L,1} \le \gamma_{L,2} \le \cdots$. The standard Riemann-von Mangoldt type formula (cf. Steuding~\cite[Theorem 7.7]{zbMATH05115272}, Smajlovi\'c~\cite{Smajlovic2010}) states that the counting function $N_L(T)$ satisfies
\[
N_L(T)=\frac{d}{2\pi} T \log T +c_L T +O(\log T),
\]
where $d>0$ is the degree of $L$ and $c_L$ is a constant depending on $L$. This matches our required asymptotic form with $C_1=d/(2\pi)$ and $C_2=\exp(-2\pi c_L/d)$. Thus, we obtain the following corollary. 
\begin{corollary}\label{cor:selberg}
Let $u$ be any non-zero real number. For a function $L \in \mathcal{S}$ with degree $d$ and constant $c_L$, if a function $f$ satisfies the conditions \textup{(C1)}--\textup{(C4)} with $C_2=\exp(-2\pi c_L/d)$, then the sequence $(uf(\gamma_{L,n}))_{\gamma_{L,n}>c}$ is uniformly distributed modulo $1$.
\end{corollary}

\subsection{Examples of functions $f$}
Having presented the relevant sequences, we now turn to explicit examples of functions $f$ that satisfy our four conditions.
Throughout this subsection, $(x_n)$ denotes any of the sequences
introduced in Subsection 2.1, that is, $(\gamma_{a,n}^{(k)})$,
$(|\gamma_{n}^{(k)}(\chi)|)$ for a primitive Dirichlet character
$\chi$, or $(\gamma_{L,n})$ for $L \in \mathcal{S}$; the constant
$C_2$ is then chosen as in the corresponding corollary in
Subsection 2.1. 
For comparison, $(q_n)$ denotes a sequence satisfying Rehberg's counting hypothesis.
\begin{example} \label{ex:15}
Let $f(t)=t^v (\log t)^w$ with $(v,w) \ne (0,0)$. The first derivative is given by
\[
f'(t)
=t^{v-1} (\log t)^{w-1} (v \log t+w).
\]
A straightforward calculation shows that for condition \textup{(C3)} to hold, we must require either $v>0$, or the case where $v=0$ and $w>1$. We thus restrict our attention to these parameter ranges.

Next, we evaluate condition \textup{(C4)}. 
Since $f(T)=T^v (\log T)^w$, the requirement $f(T)=o(T)$ is clearly satisfied if $v<1$, or if $v=1$ and $w<0$.

Combining these requirements, conditions \textup{(C3)} and \textup{(C4)} are satisfied simultaneously in the following three cases:
\begin{align}\label{uvw-condition1}
0<v<1,\quad\text{or}\quad v=0,w>1,\quad\text{or}\quad v=1,w<0.
\end{align}
Furthermore, for a sufficiently large constant $c$, the monotonicity requirements \textup{(C1)} and \textup{(C2)} are easily verified. 
Hence, for any non-zero real number $u$ and any $(v,w)$ satisfying \eqref{uvw-condition1}, the sequence $(u f(x_n))_{x_n>c}$ is uniformly distributed modulo $1$ by Theorem \ref{main}.
\end{example}
\begin{remark}
Fujii's original result for the zeros of the Riemann zeta function also covers the borderline cases $v=1, w=0$ (i.e., $f(t)=t$) and $v=0, w=1$ (i.e., $f(t)=\log t$) in Example \ref{ex:15}. In contrast, our new criterion excludes these specific cases. This exclusion is a natural consequence of the trade-off we made to avoid extensive critical-line computations and to generalize the criterion to a broader class of sequences.
\end{remark}
\begin{example}
Assume $v>1$, or $v=1$ and $w>0$. Let $f(t)=(\log t)^v (\log\log t)^w$.

Then the first derivative is given by
\[
f'(t)
=\frac{(\log t)^{v-1} (\log\log t)^{w-1} (v \log\log t+w)}{t}.
\]
For a sufficiently large constant $c$, conditions \textup{(C1)}--\textup{(C4)} are all satisfied. Therefore, for any non-zero real number $u$, the sequence $(u(\log x_n)^v (\log\log x_n)^w)$ is uniformly distributed modulo $1$. 
\end{example}
\begin{remark}
For comparison, applying Rehberg's four conditions to the sequence $(q_n)$ described in Section 1 yields the following parameter ranges. 
For $f(t)=t^v(\log t)^w$, only the case $v=0$ and $w>1$ is covered, and our criterion is strictly wider, additionally covering $0<v<1$ and the case $v=1$ and $w<0$. 
For $f(t)=(\log t)^v (\log\log t)^w$, the admissible range $v>1$, or $v=1$ and $w>0$, coincides with ours. 
This contrast reflects the difference between sequences whose counting function grows like $T\log T$ (such as $(x_n)$) and sequences whose counting function grows like $T/\log T$ (such as $(q_n)$).
\end{remark}

\section{Proofs}
In this paper, we use the following well-known classical results.
\begin{lemma}[cf.\ Lemma 4.3 of \cite{Titchmarsh1986}] \label{lem:titchmarsh_43}
Let $F(x)$ and $G(x)$ be real-valued functions, $G(x) / F^{\prime}(x)$ monotonic, and $F^{\prime}(x) / G(x) \ge m>0$, or $\le-m<0$. 
Then
\[
\left|\int_a^b G(x) e^{i F(x)} d x\right| \le \frac{4}{m}.
\]
\end{lemma}
We will prove Theorems \ref{main} and \ref{main2}. To make our proof accessible to a wider audience, we briefly recall the connection between discrete sums and Stieltjes integrals. For any sequence of real numbers $(x_n)_{n=1}^\infty$ tending to infinity with counting function $N(t,(x_n)_{n=1}^\infty)=\sum_{x_n \le t} 1$, the sum of a continuously differentiable function $g(x_n)$ can be expressed as a Riemann-Stieltjes integral: 
\[
\sum_{c < x_n \le T} g(x_n)=\int_c^T g(t) \, dN(t,(x_n)_{n=1}^\infty).
\]
Furthermore, by decomposing $N(t,(x_n)_{n=1}^\infty)$ into a smooth main term $M(t)$ and an error term $E(t)$, we can split this integral. The integral involving $dE(t)$ can then be evaluated using integration by parts for Stieltjes integrals, which essentially corresponds to Abel's summation formula.

Let $M(t)=C_1 t (\log t - \log C_2)$ be the main term of our counting function $N(t,(x_n)_{n=1}^\infty)$, and let $E(t)=O(\log t)$ be the unconditional error term. Differentiating $M(t)$ with respect to $t$, we obtain 
\[M'(t)=C_1 \left( \log t - \log C_2 +1 \right) 
=C_1 \log \left( \frac{e t}{C_2} \right).\]
Notice that this derivative precisely introduces the logarithmic factor in condition \textup{(C2)}.

By condition \textup{(C1)}, $f(t)$ is  a continuously  differentiable real-valued function on $(c,\infty)$, ensuring that the Stieltjes integral is well-defined. 
Furthermore, condition \textup{(C2)} guarantees that $f'(t) \ne 0$ for all $t>c$. Therefore, by applying Stieltjes integration to the exponential sum over the general sequence $x_n$, for a non-zero real number $h$, we have
\begin{align}
\begin{split} \label{eq:S1S2}
\sum_{c < x_n \le T} e^{i hf(x_n)}&=\int_c^T e^{i hf(t)}  dN(t,(x_n)_{n=1}^\infty) \\
&=\int_c^T e^{i hf(t)} M'(t)  dt +\int_c^T e^{i hf(t)}  dE(t) \\
&=C_1 \int_c^T e^{i hf(t)} \log \left( \frac{e t}{C_2} \right)  dt +\int_c^T e^{i hf(t)}  dE(t) \\
&\eqqcolon S_1(T) +S_2(T).
\end{split}
\end{align}

For $S_1(T)$, we apply Lemma \ref{lem:titchmarsh_43} to
\[
\int_c^T e^{i hf(t)} 
\,\log \left( \frac{e t}{C_2} \right)\,dt.
\]
In this case, it suffices that $f(t)$ and
$
\log ( e t/C_2 )
$
are real-valued, and that
\[
\frac{1}{f'(t)} \log \left( \frac{e t}{C_2} \right)
\]
is monotonic. 
By the condition (C2),
\[
\frac{f'(t)}{\log (et/C_2)}
\]
is monotonic and, on $[c,T]$, is bounded away from $0$ in the sense that
\[
\left|\frac{f'(t)}{\log (et/C_2)}\right|\ge\
\min\left\{
\left|\frac{f'(c)}{\log (ec/C_2)}
\right|,\left|\frac{f'(T)}{\log (eT/C_2)}\right|
\right\}>0.
\]
Then the required monotonicity holds and Lemma \ref{lem:titchmarsh_43} gives
\begin{align} \label{eq:int_fpp_bound}
\begin{split}
S_1(T)
&=C_1\int_{c}^{T} 
e^{i hf(t)} \log \left( \frac{e t}{C_2} \right) \,dt
=O\left(
\min\left\{
\left|\frac{hf'(c)}{\log (ec/C_2)}\right|,
\left|\frac{hf'(T)}{\log (eT/C_2)}\right|
\right\}^{-1}
\right)
\\
&=O\left(\frac{\log T}{h|f'(T)|}\right)
+O(1/h).
\end{split}
\end{align}

Now we estimate $S_2(T)$. 
We have
\begin{align}\label{uncond}
\begin{split}
S_2(T)
&=\bigl[e^{i hf(t)} E(t)\bigr]_{c}^{T}
 -\int_c^T ih\,f'(t)\,e^{i hf(t)} E(t)\,dt\\
&=e^{i hf(T)} E(T)+O(1)
 -i h\int_c^T f'(t)\,e^{i hf(t)} E(t)\,dt.
\end{split}
\end{align}
Firstly, we consider the case without assuming the Riemann Hypothesis.
By condition \textup{(C2)}, the derivative $f'(t)$ maintains a constant sign, which implies that $f(t)$ is strictly monotonic on $[c,\infty)$. If $f(t)$ were to converge to a finite limit as $t \to \infty$, the resulting sequence could not be uniformly distributed modulo $1$. Therefore, $f(t)$ must diverge to $\pm\infty$. In particular, for a sufficiently large $T$, it is guaranteed that $|f(T)|>|f(c)|$. Using this property together with $E(T)=O(\log T)$, we obtain
\begin{align}
\begin{split}
 S_2(T)
 &=O(\log T)
  +O\left(h\int_c^T |f'(t)|\,\log t \,dt\right)\\
  &=O(\log T)
  +O\left(h|f(T)|\log T\right).
  \end{split}
\label{eq:S2_final}
\end{align}

Combining \eqref{eq:int_fpp_bound} and \eqref{eq:S2_final}, we obtain
\begin{align*}
\sum_{c < x_n \le T} e^{i hf(x_n)}
&=O\left(\frac{\log T}{h|f'(T)|}\right)
+O(\log T)
+O\left(h|f(T)|\log T\right).
\end{align*}
From Weyl's criterion, it remains to show that
\begin{align*}
\frac{\sum_{c < x_n \le T} e^{i hf(x_n)}}{N(T)-N(c)} \to 0 \quad (T \to \infty).
\end{align*}
Using the earlier estimates, we have
\begin{align*}
\frac{\sum_{c < x_n \le T} e^{i hf(x_n)}}{T \log T}
&=O\left(\frac{1}{hT |f'(T)|}\right)
+O\left(\frac{1}{T}\right)
+O\left(\frac{h |f(T)| }{T }\right).
\end{align*}
Therefore, we obtain Theorem \ref{main}, since the above equality converges to $0$ under the conditions (C3) and (C4).

In the specific case where $x_n=\gamma_{0,n}^{(k)}$, by the results of Akatsuka~\cite{Akatsuka2012} for $k=1$ and Suriajaya \cite{Suriajaya2015} for $k \ge 2$, assuming the Riemann Hypothesis yields the improved error bound $E(T)=O\left( \frac{\log T}{\sqrt{\log\log T}} \right)$. By substituting this bound to \eqref{uncond} and following the same reasoning as above, we obtain Theorem \ref{main2}.

We now give the proof of Theorem \ref{thm:discrepancy}. 
By applying the Erd\H{o}s--Tur\'an inequality together with our earlier estimates, we obtain
\begin{align*}
D_{N(T)-N(c)}
&=O\left(\frac{1}{H}\right)
+O\left(\frac{1}{T |f'(T)|}\right)
+O\left(\frac{\log H}{T}\right)+O\left(\frac{H|f(T)| }{T }\right).
\end{align*}
Therefore, by choosing
\[
H=\left\lfloor \left( \frac{T}{|f(T)| } \right)^{1/2} \right\rfloor,
\]
we arrive at the desired bound. This completes the proof of Theorem \ref{thm:discrepancy}.

\section*{Acknowledgments}
The authors would like to express their sincere gratitude to the anonymous referee for the careful reading of the manuscript and for the many valuable comments and suggestions, which have significantly improved the quality and presentation of this paper.
This work was supported by JSPS KAKENHI Grant Number JP22K13897 (Murahara).

\bibliographystyle{amsalpha}
\bibliography{References} 
\end{document}